\documentclass{amsart}
\usepackage{appendix}
\usepackage{graphicx}
\usepackage[parfill]{parskip}
\usepackage{amsfonts, amscd, mathrsfs, amsmath, amssymb, amsthm}
\usepackage{url}
\usepackage{tikz}
\usepackage{color, verbatim}
\usepackage{bm}
\usepackage[hmargin=3cm,vmargin=3cm]{geometry}
\usepackage{hyperref}
\hypersetup{pdfpagemode=UseNone,colorlinks=true}
\numberwithin{equation}{section}

\newtheorem{theorem}[equation]{Theorem}

\newtheorem{lemma}[equation]{Lemma}

\newtheorem{definition}[equation]{Definition}
\theoremstyle{definition}

\theoremstyle{remark}

\DeclareMathOperator {\Diff} {Diff}
\DeclareMathOperator {\Met} {Met}

\DeclareMathOperator {\Neg} {Neg}

\DeclareMathOperator {\area} {area}

\begin{document}

\title{On moduli spaces of uniformly negatively curved metrics}
\author{Yasha Savelyev}
\thanks {Supported by CONAHCYT research grant CF-2023-I}
\email{yasha.savelyev@gmail.com}
\address{University of Colima}

\begin{abstract}
Let $\Neg(X)$ denote the space of complete Riemannian metrics
with uniformly negative curvature on a surface $X$, equipped
with the intrinsic uniform $C^2$ topology. Let
$\mathcal M(X)=\Neg(X)/\Diff(X)$ be the corresponding moduli space,
with the quotient topology. We construct elementary locally
constant functionals
on $\mathcal M(X)$, with values in finite symmetric products of
$\{0,1\}$, based on geodesic string counts. As an upshot we show that
$\mathcal M(\mathbb{R}\times S^1)$ is disconnected. This is perhaps surprising: the two metrics we separate
are joined by an explicit path of metrics with constant curvature
$-1$. The point is that this path is only continuous in the weak
Whitney topology. More generally, if $X$ is a finite-type surface of
hyperbolic type with $n$ punctures, then the pure moduli space has at
least $2^n$ connected components, while $\mathcal M(X)$ has at least
$n+1$ connected components. 
\end{abstract}
\maketitle

\section{Introduction}
Spaces and moduli spaces of
negatively curved metrics have been studied before, especially for
closed manifolds, see for instance
~\cite{cite_FarrellOntanedaNegativeMetrics},
~\cite{cite_FarrellOntanedaTeich}. For the analogous questions for
spaces of metrics with positive scalar curvature and with non-negative
curvature, see ~\cite{cite_SpaceOfPositiveCurvaturemetrics} and
~\cite{cite_tuschmann2021moduli}, respectively. The noncompact surface case with
the intrinsic uniform \(C^2\) topology seems to be
a different and less explored setting. Our main thrust
is to prove disconnection results for the moduli spaces of
metrics on surfaces. Our approach is
exceptionally simple: using Gauss--Bonnet and Arzelà-Ascoli
theorems to
establish compactness of certain closed geodesic spaces, and
constructing locally constant geodesic counting functionals.  At the same time the idea shadows
a deeper theory involving the Fuller index, which we will indicate.

Let $\Neg(X)$ denote the space of Riemannian, complete
metrics $g$ with uniformly negative curvature on a surface $X$. By
uniformly negative we mean that there is a constant $a>0$ such that
the Gaussian curvature satisfies $K_g\leq -a$ on all of $X$.
The space $\Neg(X)$ will be
given the intrinsic uniform $C^2$ topology, see Definition
\ref{def_uniform}. This is a metrizable topology, finer than the
weak Whitney $C^2_{\mathrm{loc}}$ topology but coarser than
the strong Whitney $C ^{2}$ topology. The corresponding moduli space is
\begin{equation*}
\mathcal M(X)=\Neg(X)/\Diff(X),
\end{equation*}
with the quotient topology, where diffeomorphisms act by pullback. If
$X$ has distinguished punctures, we also write
\begin{equation*}
\mathcal M_{\mathrm p}(X)=\Neg(X)/\Diff_{\mathrm p}(X),
\end{equation*}
where $\Diff_{\mathrm p}(X)$ denotes the subgroup preserving each
puncture, or equivalently each end, individually.

This moduli space should be distinguished from the usual Riemann
moduli space $\mathcal M_{g,n}$. By uniformization, the latter may be
viewed as the moduli space of complete finite-area hyperbolic metrics,
so all punctures are cusp ends. In the present paper the ambient moduli
space is larger: the same underlying smooth punctured surface is
allowed to carry uniformly negatively curved metrics with either cusp
or funnel type behavior at a puncture. The intrinsic uniform topology
remembers this large-scale end geometry. Thus the all-cusp hyperbolic
locus lies in one sector of our space, while the geodesic
string counting invariant below
separates the other cusp/funnel sectors. 
This should be compared with the usual Teichm\"uller topology on the
all-cusp hyperbolic locus.
\begin{theorem}\label{thm_cylinder_moduli}
The moduli space
\begin{equation*}
\mathcal M(\mathbb{R}\times S^1)
=
\Neg(\mathbb{R}\times S^1)/\Diff(\mathbb{R}\times S^1)
\end{equation*}
is disconnected.
\end{theorem}

\paragraph{The role of the topology} \label{rem_weakpath}
The disconnection in Theorem \ref{thm_cylinder_moduli} is genuinely a
feature of the intrinsic uniform topology. The two metrics $g_0,g_1$ exhibited
in Section \ref{sec_Acalculation}, with $N(g_0,\beta)=0$ and
$N(g_1,\beta)=1$, are the warped products
\begin{equation*}
g_i=dx^2+f_i(x)^2\,d\theta^2,
\end{equation*}
with
\begin{equation*}
f_0(x)=e^x,\qquad f_1(x)=\cosh x.
\end{equation*}
Interpolating the warping functions linearly,
\begin{equation*}
f_s(x)=(1-s)e^x+s\cosh x,\qquad s\in[0,1],
\end{equation*}
one has $f_s''=f_s$, and therefore
\begin{equation*}
K_{g_s}=-\frac{f_s''}{f_s}=-1.
\end{equation*}
Thus each $g_s=dx^2+f_s(x)^2d\theta^2$ lies in
$\Neg(\mathbb R\times S^1)$. The path $s\mapsto g_s$ is
continuous in the weak Whitney ($C^\infty_{\mathrm{loc}}$) topology,
so the two metrics lie in one weak Whitney path component; in
particular their images in $\mathcal M(\mathbb R\times S^1)$ lie
in one path component for the quotient of the weak Whitney topology. 

At the marked level the particular metrics in this example are also
not globally bi-Lipschitz, and so the order-zero part of the topology
already separates them. For this reason the applications below are
stated on the moduli spaces. A diffeomorphism of a non-compact surface
need not preserve a global bi-Lipschitz class, while the invariant
constructed here is natural under diffeomorphisms.

We now describe a basic metric deformation invariant.
A \emph{geodesic string} in
$X,g$ will be short for: an equivalence class of a constant speed
smooth $g$-geodesic mapping $S^1\to X$, under the $S^1$
reparametrization action. Let $\pi^*_1(X)$ be the set of
non-trivial free homotopy classes of loops in $X$
(topologized with the
discreet topology).
Define
\begin{equation*}
N: \Neg(X)\times \pi^*_1(X)\to \{0,1\},
\end{equation*}
Define \[
N:\Neg(X)\times \pi^*_1(X)\to \{0,1\}
\]
by declaring \(N(g,\beta)\) to be the number of class \(\beta\)
\(g\)-geodesic strings. This number is at most one by the negative
curvature assumption. 

\begin{theorem}\label{thm_mainN}
The function $N$ is locally constant. It is natural under diffeomorphisms: if
$\phi:X\to X$ is a diffeomorphism, then
\begin{equation*}
N(\phi^*g,\beta)=N(g,\phi_*\beta).
\end{equation*}
Finally,
\begin{equation*}
N(g,\beta)=N(g,\beta^{-1}).
\end{equation*}
Consequently, for every subgroup $G\subset \Diff(X)$, $N$ descends to
a locally constant function
\begin{equation*}
N_G:(\Neg(X)\times \pi^*_1(X))/G\to \{0,1\}.
\end{equation*}
Here $G$ acts diagonally by
\begin{equation*}
(g,\beta)\cdot \phi=(\phi^*g,(\phi^{-1})_*\beta).
\end{equation*}
In particular, if $\mathcal B\subset \pi^*_1(X)$ is a finite subset
which is preserved by $G$ up to permutation and inversion, then the
unordered collection of values defines a locally constant function
\begin{equation*}
\mathcal N_{\mathcal B,G}:\Neg(X)/G\to
\operatorname{Sym}^{|\mathcal B|}(\{0,1\}),
\qquad
\mathcal N_{\mathcal B,G}([g])=
\{N(g,\beta)\}_{\beta\in\mathcal B}.
\end{equation*}
If $|\mathcal B|=1$, we regard $\mathcal N_{\mathcal B,G}$ as a
$\{0,1\}$-valued function.
\end{theorem}
\paragraph{A deeper theory in the background}
Let $\Met (X)$ denote the space of complete Riemannian
metrics with uniformly negatively curved ends.
In forthcoming work, we construct an extension of $N$ to
a locally constant, $\Diff (X)$ natural functional:
$$F:\Met(X)\times \pi^*_1(X) \to \mathbb{Q},$$ built using
Fuller index theory, but still based on the same idea of
counting geodesic strings. 
This requires a more involved analogue of the
compactness Theorem \ref{thm_compactness}, and is beyond the
methods of the present paper. Granting it, we would get that $\Met(\mathbb{R}\times
S^1)/\Diff(\mathbb{R}\times S^1) $ is disconnected.

Going back to the theorem above, taking $\mathcal B$ to be the
peripheral classes around the ends, the local constancy
of the functions
$\mathcal N_{\mathcal B,G}$ gives the following
finite-type moduli statement.
\begin{theorem}\label{thm_finite_type_moduli}
Let $X$ be an orientable finite-type surface without boundary, with
$n\geq1$ punctures, and assume that $X$ is of hyperbolic type, i.e.
$\chi(X)<0$. Then the pure moduli space
\begin{equation*}
\mathcal M_{\mathrm p}(X)=\Neg(X)/\Diff_{\mathrm p}(X)
\end{equation*}
has at least $2^n$ connected components, and the full moduli space
\begin{equation*}
\mathcal M(X)=\Neg(X)/\Diff(X)
\end{equation*}
has at least $n+1$ connected components. 
\end{theorem}

\section{\texorpdfstring{The space of geodesic strings and $N$}{The space of geodesic strings and N}}
\label{sec_geodesic_strings}
Let $X$ be a surface admitting a metric in $\Neg (X)$, and
fix a non-trivial class $\beta \in \pi ^{*} _{1} (X)$. Free homotopy
classes of loops correspond to conjugacy classes of deck
transformations of the universal cover $\widetilde X \to X$; let
$\eta _{\beta}$ be a deck transformation in the conjugacy class
determined by $\beta$. It is well defined up to conjugacy, and
non-triviality of $\beta$ means $\eta _{\beta} \neq \mathrm{id}$.
The translation length
\begin{equation*}
|\eta _{\beta}| _{g}=\inf _{\widetilde x\in\widetilde X}
\widetilde d(\widetilde x,\eta _{\beta}\widetilde x)
\end{equation*}
used below is conjugation invariant, so depends only on $\beta$. In
the cylinder application of Sections \ref{sec_Acalculation}--
\ref{sec_cylinder_moduli} we specialise to $X=\mathbb{R}\times S^1$, with
$\beta$ the generator of $\pi_1(X)\simeq\mathbb{Z}$. In what follows
all metrics are elements of $\Neg(X)$, with the latter as in
the Introduction.

Let $LX$ denote the space of all smooth maps
$S^1\to X$,
\begin{equation*}
S^1=\mathbb{R}/\mathbb{Z},
\end{equation*} 
with the usual $C^0$ topology. The
space $S(g,\beta)\subset LX$ will denote the subspace of all
constant speed parametrized, $g$-geodesics in class $\beta$. Set
\begin{equation*}
\mathcal{O}(g,\beta)=S(g,\beta)/S^1,
\end{equation*}
the quotient by the $S^1$ action, where $S^1$ is acting by
reparametrization: $t\cdot o(\tau)=o(\tau+t)$. The elements of
$\mathcal{O}(g,\beta)$ are called geodesic strings.

Since the metrics in $\Neg(X)$ have strictly negative curvature, there
is at most one class $\beta$ geodesic string. Indeed, two geometrically
distinct such geodesics would lift to axes of the same deck
transformation in the Hadamard universal cover, and the flat strip
theorem would give a flat strip, impossible under negative curvature.
We therefore define
\begin{equation*}
N(g,\beta)=\#\mathcal O(g,\beta)\in\{0,1\}.
\end{equation*}
Theorem \ref{thm_mainN} says that this elementary count is a locally constant deformation invariant.

We give $\Neg(X)$ the intrinsic uniform $C^2$ topology,
which is very natural in our setup, and particularly in the
compactness argument of Theorem \ref{thm_compactness}.
\begin{definition}\label{def_uniform}
The intrinsic uniform \(C^2\) topology is defined as follows. A
neighbourhood basis of \(g\in\Neg(X)\) consists,
for $\lambda>1$ and $\epsilon>0$, of the sets
\begin{equation*}
N_{g,\lambda,\epsilon}=\Big\{ h\in\Neg(X)\;\Big|\;
\sup_{x\in X}\sup_{0\neq v\in T_xX}
\Big|\log\frac{h_x(v,v)}{g_x(v,v)}\Big|<\log\lambda,
\;
\sup_{x\in X}\big(|\nabla^g(h-g)|_g+|(\nabla^g)^2(h-g)|_g\big)
<\epsilon\Big\},
\end{equation*}
where $\nabla^g$, $|\cdot|_g$ are the Levi-Civita connection and
tensor norms of $g$. That is, $h$ is uniformly $\lambda$-bi-Lipschitz
to $g$ at order $0$, and uniformly $\epsilon$-close at orders $1$ and
$2$, measured in $g$ itself. This is a metrizable topology, stronger
than the weak Whitney $C^2_{\mathrm{loc}}$ topology. The order-$0$ condition is
multiplicative on purpose: it controls the length functional
uniformly over all of $X$, including along curves whose images leave
every compact set, and is the only order used in Lemma
\ref{lem_Lcontinuous}. Order $1$ controls the Christoffel symbols,
hence which curves are nearly geodesic, see Lemma \ref{lem_kappa}.
Order $2$ controls the curvature, hence the geodesic flow and its
linearisation, which is what the persistence of closed geodesics in
Theorem \ref{thm_mainN} uses. No derivative of order $\geq 3$ is
used anywhere.
\end{definition}

For a fixed metric the compactness we are after is elementary, and it
is instructive to isolate why: a class $\beta$ geodesic is
automatically a minimiser, so it can neither collapse nor, under the uniform negative curvature assumption, spread over a
non-compact region. Theorem \ref{thm_compactness} is the assertion
that this confinement survives passage to a compact family of
metrics.

\begin{lemma}\label{lem_confinement}
Fix $g\in\Neg(X)$. Then every class $\beta$ $g$-geodesic
string has length equal to the translation length
\begin{equation*}
|\eta_\beta|_g:=\inf_{\widetilde x\in\widetilde X}
\widetilde d(\widetilde x,\eta_\beta\widetilde x),
\end{equation*}
and is a minimiser of length in $\beta$; in particular all such
geodesics have the same length, and
\begin{equation*}
\inf_{o\in\beta}\ell_g(o)=|\eta_\beta|_g.
\end{equation*}
Moreover the images of all elements of $S(g,\beta)$ are contained in
a single compact subset of $X$.
\end{lemma}

\begin{proof}
Let $\pi:\widetilde X\to X$ be the universal cover, a Hadamard
manifold by Cartan--Hadamard, and regard $\eta_\beta$ as a deck
isometry. That a closed geodesic in a non-positively curved manifold
minimises length in its free homotopy class, with length equal to the
translation length $|\eta_\beta|_g$ of the corresponding deck
isometry, is well known; it follows from the standard correspondence
between such geodesics and axes of $\eta_\beta$ and the convexity of
the displacement function on a Hadamard manifold, see for instance
~\cite[II.6]{cite_BridsonHaefliger}. In particular all class $\beta$
geodesics have the common length $|\eta_\beta|_g$ and
$\inf_{o\in\beta}\ell_g(o)=|\eta_\beta|_g$.

For the confinement, there is nothing to prove if $S(g,\beta)$ is
empty. If there is only one geometric class $\beta$ geodesic, its
image is compact. Otherwise choose two geometrically distinct class
$\beta$ geodesics $\gamma_0,\gamma_1$, with lifts
$\widetilde\gamma_0,\widetilde\gamma_1$. By the preceding paragraph,
these axes have equal translation length, hence are parallel, and by
the flat strip theorem
~\cite[Proposition 5.1]{cite_FlatStripEberleinOneil} they bound a
flat, totally geodesic strip in $\widetilde X$. This is impossible
because $K_g\leq -a<0$ on all of $X$, for some $a>0$. Hence there is
at most one class $\beta$ geodesic string, and the images of all
elements of $S(g,\beta)$ are contained in a single compact subset of
$X$.
\end{proof}

Given a subset $\mathscr{K}\subset\Neg(X)$ set
\begin{equation*}
\mathcal{O}(\mathscr{K},\beta):=\{(o,g)\in L X\times\mathscr{K}\mid
o\in \mathcal{O} (g,\beta)\}.
\end{equation*}

\begin{theorem}\label{thm_compactness}
Let $\mathscr{K}\subset\Neg(X)$ be compact. Then the space
$\mathcal{O} (\mathscr{K},\beta)$ is compact.
\end{theorem}

The proof uses a few lemmas. The first is the continuity of the
marginal length; it is exactly here that the uniform, rather than
weak Whitney, topology is used.

\begin{lemma}\label{lem_Lcontinuous}
The marginal length
\begin{equation}\label{eq_bound}
\mathcal{L}(g)=\inf_{o\in\beta}\ell_g(o)
\end{equation}
defines a finite continuous function on $\Neg(X)$, where
\begin{equation*}
\ell_g(o)=\int_{S^1}\sqrt{g(\dot o,\dot o)}.
\end{equation*}
\end{lemma}

\begin{proof}
Fixing any smooth $o_\beta\in\beta$ gives
$\mathcal{L}(g)\leq\ell_g(o_\beta)<\infty$. For the continuity, fix
$g$. By Definition \ref{def_uniform}, for every $\delta>0$ there is
a neighbourhood $V$ of $g$ such that, for all $h\in V$,
\begin{equation}\label{eq_bilip}
(1+\delta)^{-1}g(v,v)\leq h(v,v)\leq (1+\delta)g(v,v)
\qquad \text{for all }x\in X,\;v\in T_xX.
\end{equation}
Taking square roots and integrating along any loop $o$ gives
\begin{equation*}
(1+\delta)^{-1}\ell_g(o)\leq\ell_h(o)\leq
(1+\delta)\ell_g(o),
\end{equation*}
with a bound independent of $o$. Since this estimate holds for every
$o\in\beta$, it passes to the infimum:
\begin{equation*}
(1+\delta)^{-1}\mathcal{L}(g)\leq\mathcal{L}(h)\leq
(1+\delta)\mathcal{L}(g),\qquad h\in V.
\end{equation*}
Hence $\mathcal{L}$ is continuous at $g$. Note that no closed
geodesic need exist for this argument: only near-minimisers are used,
so it applies verbatim where the infimum \eqref{eq_bound} is not
attained.
\end{proof}

The next lemma records the only fact about the limit metric needed
below. It gives a positive lower length bound, but does not assert
that this lower bound is attained.

\begin{lemma}\label{lem_lower_length}
Let $g_j\to g'$ in $\Neg(X)$, and let
$b_j\in S(g_j,\beta)$. Then
\begin{equation*}
\mathcal{L}(g')=\inf_{o\in\beta}\ell_{g'}(o)>0.
\end{equation*}
\end{lemma}

\begin{proof}
By Lemma \ref{lem_confinement}, $b_j$ minimises $g_j$-length in the
class $\beta$. Hence, for every loop $o\in\beta$,
\begin{equation*}
\mathcal{L}(g_j)=\ell_{g_j}(b_j)
\leq \ell_{g_j}(o).
\end{equation*}
For $j$ large, $g_j$ and $g'$ are uniformly bi-Lipschitz close, so
\begin{equation*}
\ell_{g_j}(o)\leq (1+\delta)\ell_{g'}(o)
\end{equation*}
for all $o\in\beta$. Taking the infimum over $o$ gives
\begin{equation*}
\mathcal{L}(g_j)\leq (1+\delta)\mathcal{L}(g').
\end{equation*}
Since $b_j$ is a genuine closed geodesic, $\mathcal{L}(g_j)>0$. If
$\mathcal{L}(g')$ were zero, the preceding inequality would force
$\mathcal{L}(g_j)=0$ for all large $j$, a contradiction. Hence
$\mathcal{L}(g')>0$.
\end{proof}

The second lemma quantifies the elementary fact that a geodesic of a
nearby metric is nearly geodesic; this is where order $1$ of
Definition \ref{def_uniform} is used.

\begin{lemma}\label{lem_kappa}
Let $g,h\in\Neg(X)$ with $h\in N_{g,\lambda,\epsilon}$ as in
Definition \ref{def_uniform}, and let $b$ be a closed $h$-geodesic.
Then the $g$-geodesic curvature of $b$ satisfies
\begin{equation*}
|\kappa_g(b)|\leq C(\lambda)\epsilon,
\end{equation*}
where $C(\lambda)$ depends only on $\lambda$.
\end{lemma}

\begin{proof}
The difference of Levi-Civita connections
$\Lambda:=\nabla^g-\nabla^h$ is a tensor, given by the standard
formula
\begin{equation*}
2h(\Lambda(u,v),w)=-(\nabla^g_uh)(v,w)- (\nabla^g_vh)(u,w)
+(\nabla^g_wh)(u,v),
\end{equation*}
and $\nabla^g h=\nabla^g(h-g)$. Using the bi-Lipschitz comparison
between $|\cdot|_g$ and $|\cdot|_h$ this gives
\begin{equation*}
|\Lambda|_g\leq C(\lambda)\sup_X|\nabla^g(h-g)|_g
\leq C(\lambda)\epsilon.
\end{equation*}
If $b$ is an $h$-geodesic then
$\nabla^g_{\dot b}\dot b=\Lambda(\dot b,\dot b)$, so
\begin{equation*}
|\kappa_g(b)|\leq \frac{|\nabla^g_{\dot b}\dot b|_g}{|\dot b|_g^2}
\leq |\Lambda|_g\leq C(\lambda)\epsilon.
\end{equation*}
\end{proof}

The third lemma is the geometric heart of the argument, an
approximate form of the flat strip theorem. It is proved by passing
to the cyclic cover of $X$ associated to the subgroup
$\langle\eta_\beta\rangle\subset\pi_1(X)$, and using Gauss--Bonnet.

\begin{lemma}\label{lem_nearlyflat}
Let $(X,g)$ be a complete surface with $K_g\leq0$, and let
$b_0,b_1:S^1\to X$ be disjoint loops in the same non-trivial free
homotopy class $\beta $. Suppose that, for $i=0,1$, the loop $b_i$ has
$g$-length $<D$, $g$-geodesic curvature $<\delta_i$, and an
$\eta _{\beta}$-equivariant lift to the universal cover which is an embedded
line. Let $\widehat\Omega$ be the compact annulus in the cyclic cover
cobounded by the corresponding lifts $\widehat b_0,\widehat b_1$.
Let $U\subset\widehat\Omega$ be an open set separating
$\widehat b_0,\widehat b_1$, disjoint from both, of positive finite
area. Then there is a point $x\in U$ with
\begin{equation*}
K_g(x)>-\epsilon(U,D,\delta_0,\delta_1),
\end{equation*}
where
\begin{equation*}
\lim_{\delta_0,\delta_1\to0}\epsilon(U,D,\delta_0,\delta_1)=0.
\end{equation*}
Here $K_g$ denotes the curvature of the lifted metric on the cyclic
cover.
\end{lemma}
\begin{proof}
Let
$p:\widehat X\to X$ be the
quotient of the universal cover by the action of $\langle
\eta _{\beta} \rangle $. With the pulled-back metric,
$p$ is a local isometry and $\widehat X$ is an infinite cylinder. The
loops $b_0,b_1$ lift to closed homotopic loops
$\widehat b_0,\widehat b_1$ in $\widehat X$.

By hypothesis $\eta _{\beta }$-equivariant  lifts of $b_0,b_1$ to the universal
cover are embedded. Hence $\widehat b_0$ and $\widehat b_1$ are
embedded circles in $\widehat X$. Since the original loops are
disjoint, after choosing the relevant lifts these two circles are
disjoint and cobound a compact sub-annulus $\widehat\Omega$.

\emph{Gauss--Bonnet on the cobounding annulus.} By Gauss--Bonnet,
\begin{equation*}
\int_{\widehat\Omega}K_g\,dA+
\int_{\partial\widehat\Omega}\kappa_g\,ds
=2\pi\chi(\widehat\Omega)=0.
\end{equation*}
The boundary curves have lengths $<D$ and geodesic curvatures
$<\delta_0$, $<\delta_1$, as $p$ is a local isometry. Therefore
\begin{equation*}
\left|\int_{\partial\widehat\Omega}\kappa_g\,ds\right|
<D(\delta_0+\delta_1),
\end{equation*}
and consequently
\begin{equation}\label{eq_GBtwoboundary}
\left|\int_{\widehat\Omega}K_g\,dA\right|<D(\delta_0+\delta_1).
\end{equation}
Since $K_g\leq0$ on $\widehat\Omega$ and $U\subset\widehat\Omega$,
\begin{equation*}
0\leq -\int_UK_g\,dA
\leq -\int_{\widehat\Omega}K_g\,dA
<D(\delta_0+\delta_1).
\end{equation*}
If
\begin{equation*}
K_g(x)\leq -\frac{D(\delta_0+\delta_1)}{\area(U)}
\qquad\text{for every }x\in U,
\end{equation*}
then
\begin{equation*}
-\int_UK_g\,dA\geq D(\delta_0+\delta_1),
\end{equation*}
contradicting the previous strict inequality. Hence there is a point
$x\in U$ such that
\begin{equation*}
K_g(x)>-\frac{D(\delta_0+\delta_1)}{\area(U)}.
\end{equation*}
Taking
\begin{equation*}
\epsilon(U,D,\delta_0,\delta_1)=
\frac{D(\delta_0+\delta_1)}{\area(U)},
\end{equation*}
we are done.
\end{proof}

\begin{proof}[Proof of Theorem \ref{thm_compactness}]
Since the topology on $\Neg(X)$ is metrizable, it suffices to prove
sequential compactness of
\begin{equation*}
S(\mathscr{K},\beta):=\{(o,g)\in L X\times\mathscr{K}\mid
o\in S (g,\beta)\}.
\end{equation*}
Let
\begin{equation*}
(b_j,g_j)\in S (\mathscr{K},\beta)
\end{equation*}
be a sequence. Passing to a subsequence, we may assume that
$g_j\to g'\in\mathscr{K}$. We show that the corresponding loops have
a $C^0$-convergent subsequence.

First, the lengths are controlled with respect to the fixed metric
$g'$. Choose a smooth loop $o_\beta\in\beta$. Since $g_j\to g'$, the
metrics $g_j$ and $g'$ are uniformly bi-Lipschitz for all large $j$, and
\begin{equation*}
\ell_{g_j}(b_j)=\mathcal L(g_j)
\leq \ell_{g_j}(o_\beta)
\end{equation*}
is bounded. Enlarging the constant, we therefore have
\begin{equation}\label{eq_uniform_D}
\ell_{g'}(b_j)<D
\end{equation}
for all large $j$. Since the $b_j$ are parametrised with constant
$g_j$-speed and $g_j$ is uniformly bi-Lipschitz to $g'$, the family is
equicontinuous with respect to the distance of $g'$.

If the images of the $b_j$ are contained in a fixed compact subset of
$X$, Arzel\`a--Ascoli gives a $C^0$-convergent subsequence. The
geodesic equation, together with the uniform $C^2$ convergence
$g_j\to g'$ on this compact set, implies that the limit is a smooth
$g'$-geodesic in the same free homotopy class. Thus it remains only to
rule out escape to infinity.

Suppose, for contradiction, that after passing to a subsequence the
images of the $b_j$ leave every compact subset of $X$. By Lemma
\ref{lem_kappa} and the convergence $g_j\to g'$,
\begin{equation}\label{eq_kappadecay}
\delta_j:=|\kappa_{g'}(b_j)|_{C^0}\to0.
\end{equation}
By Lemma \ref{lem_lower_length}, $\mathcal L(g')>0$. Thus every loop
in the class $\beta$ has $g'$-length at least $\mathcal L(g')$.
Since $g'\in\Neg(X)$, there is $a>0$ such that
\begin{equation*}
K_{g'}\leq -a
\end{equation*}
on all of $X$.

The $g'$-diameters of the curves $b_j$ are bounded by $D$. Since their
images leave every compact set, after passing to a subsequence we may
choose indices $j<k$ so that
\begin{equation*}
d_{g'}(b_j,b_k)>2.
\end{equation*}
We also choose $j,k$ so large that $\delta_j$ and $\delta_k$ are small.

Let $p:\widehat X\to X$ be the cyclic cover associated to
$\langle\eta_\beta\rangle$, and let $\widehat\Omega$ be the compact
annulus cobounded by the corresponding lifts of $b_j$ and $b_k$. Let
$r$ be the $\widehat g'$-distance from the lift of $b_j$, restricted
to $\widehat\Omega$, and set
\begin{equation*}
U=\widehat\Omega\cap r^{-1}((0,1)).
\end{equation*}
Then $U$ is an open set separating the two boundary components and is
disjoint from both. By the coarea formula, for almost every
$t\in(0,1)$ the level $r^{-1}(t)\cap\widehat\Omega$ has a component
separating the two boundary circles. 
Such a separating component represents the generator of the cyclic cover, and hence projects to a loop in the class \(\beta\) or \(\beta^{-1}\). Since these two classes have the same marginal length,
its length is at least \(\mathcal L(g')\).
Thus
\begin{equation*}
\area(U)\geq \mathcal L(g').
\end{equation*}

The $\eta_\beta$-equivariant lifts of $b_j$ and $b_k$ are embedded
lines, since they are geodesics for the Hadamard metrics
$\widetilde g_j$ and $\widetilde g_k$. Thus the hypotheses of Lemma
\ref{lem_nearlyflat}, applied to the metric $g'$, the loops $b_j,b_k$,
the lifted annulus and this open set $U$, are satisfied. Since
\begin{equation*}
\epsilon(U,D,\delta_j,\delta_k)
\leq \frac{D(\delta_j+\delta_k)}{\mathcal L(g')},
\end{equation*}
we may take $j,k$ so large that this number is $<a$. The lemma gives a
point $x\in U$ with
\begin{equation*}
K_{g'}(x)>-a,
\end{equation*}
contradicting $K_{g'}\leq -a$ on all of $X$. This rules out escape to
infinity, and hence proves compactness.
\end{proof}

\begin{proof}[Proof of Theorem \ref{thm_mainN}]
We now prove local constancy. Fix $\beta\in\pi_1^*(X)$ and
$g\in\Neg(X)$. 
Suppose first that $N(g,\beta)=0$. If $N(\cdot,\beta)$ were not
locally zero near $g$, then there would be metrics $g_j\to g$ and
class $\beta$ geodesic strings $o_j\in\mathcal O(g_j,\beta)$. The set
\begin{equation*}
\{g\}\cup\{g_j\}_{j\geq1}
\end{equation*}
is compact in $\Neg(X)$. By Theorem \ref{thm_compactness}, after
passing to a subsequence the strings $o_j$ converge to a class $\beta$
$g$-geodesic string, contradicting $N(g,\beta)=0$. Thus
$N(\cdot,\beta)$ is locally zero near $g$.

Now suppose that $N(g,\beta)=1$, and let \(o\) be the unique class \(\beta\) geodesic string. Since
\(K_g<0\), it has no non-trivial periodic normal Jacobi
field, and is therefore non-degenerate, as is well known.
By the implicit function theorem for non-degenerate
closed geodesics, the geodesic string $o$ persists uniquely for all
metrics sufficiently close to $g$. Thus $N(\cdot,\beta)=1$ on a
neighbourhood of $g$.

This proves that
\begin{equation*}
N:\Neg(X)\times\pi_1^*(X)\to\{0,1\}
\end{equation*}
is locally constant.

Naturality is immediate. A diffeomorphism $\phi:X\to X$ identifies the
$\phi^*g$-geodesic strings in the class $\beta$ with the $g$-geodesic
strings in the class $\phi_*\beta$. Hence
\begin{equation*}
N(\phi^*g,\beta)=N(g,\phi_*\beta).
\end{equation*}
Reversing the orientation of a closed geodesic identifies the class
$\beta$ strings with the class $\beta^{-1}$ strings, so
\begin{equation*}
N(g,\beta)=N(g,\beta^{-1}).
\end{equation*}

For the diagonal action in the statement,
\begin{equation*}
(g,\beta)\cdot\phi=(\phi^*g,(\phi^{-1})_*\beta),
\end{equation*}
the naturality formula gives
\begin{equation*}
N(\phi^*g,(\phi^{-1})_*\beta)=N(g,\beta).
\end{equation*}
Therefore $N$ descends to
\begin{equation*}
N_G:(\Neg(X)\times\pi_1^*(X))/G\to\{0,1\}.
\end{equation*}
The quotient map is open, since it is the quotient by a group action by
homeomorphisms. Thus the descent of a locally constant invariant is
again locally constant.

The assertion about the functions $\mathcal N_{\mathcal B,G}$ is the
same statement applied to a finite collection of classes, using
$N(g,\beta)=N(g,\beta^{-1})$ when a diffeomorphism reverses a class.
\end{proof}

\section{Sample calculations} \label{sec_Acalculation}
Let $X=\mathbb{R}\times S^1$, and let $\beta$ be the class of the
loop $\gamma(t)=(0,t)$. For a positive function $f$ on $\mathbb{R}$,
consider the warped product metric
\begin{equation*}
g_f=dx^2+f(x)^2d\theta^2.
\end{equation*}
Here $d\theta^2$ is the standard metric on $S^1=\mathbb{R}/\mathbb{Z}$;
for background on warped products, see
~\cite{cite_ManifoldsOfNegativeCurvatureBishopOneil}. For such a
metric the Gauss curvature is
\begin{equation}\label{eq_warpcurvature}
K_{g_f}=-\frac{f''}{f}.
\end{equation}
Also, the circle $x=x_0$ is a geodesic if and only if
\begin{equation}\label{eq_circlegeodesic}
f'(x_0)=0.
\end{equation}
Indeed, the only relevant Christoffel symbol is
\begin{equation*}
\Gamma^x_{\theta\theta}=-ff',
\end{equation*}
so a circle $x=x_0$ is geodesic exactly when $f'(x_0)=0$. Every loop
in the generator class is freely homotopic to one of these circles,
and the length of the circle $x=x_0$ is proportional to $f(x_0)$.

Now set
\begin{equation*}
f_0(x)=e^x,\qquad f_1(x)=\cosh x,
\end{equation*}
and
\begin{equation*}
g_0=dx^2+e^{2x}d\theta^2,
\qquad
g_1=dx^2+\cosh^2(x)d\theta^2.
\end{equation*}
Since $f_0''=f_0$ and $f_1''=f_1$, equation
\eqref{eq_warpcurvature} gives
\begin{equation*}
K_{g_0}=K_{g_1}=-1.
\end{equation*}
Thus both metrics belong to $\Neg(X)$.

For $g_0$, one has $f_0'(x)=e^x>0$ for all $x$. Hence there is no
closed geodesic in the generator class, and so
\begin{equation*}
N(g_0,\beta)=0.
\end{equation*}
For $g_1$, one has $f_1'(x)=\sinh x$, so the unique circle satisfying
\eqref{eq_circlegeodesic} is $x=0$. Thus $g_1$ has a single class
$\beta$ geodesic string, represented by $t\mapsto(0,t)$. Therefore
\begin{equation*}
N(g_1,\beta)=1.
\end{equation*}

\begin{figure}[h]
\centering
\begin{minipage}[t]{0.48\textwidth}
\centering
\includegraphics[width=\textwidth]{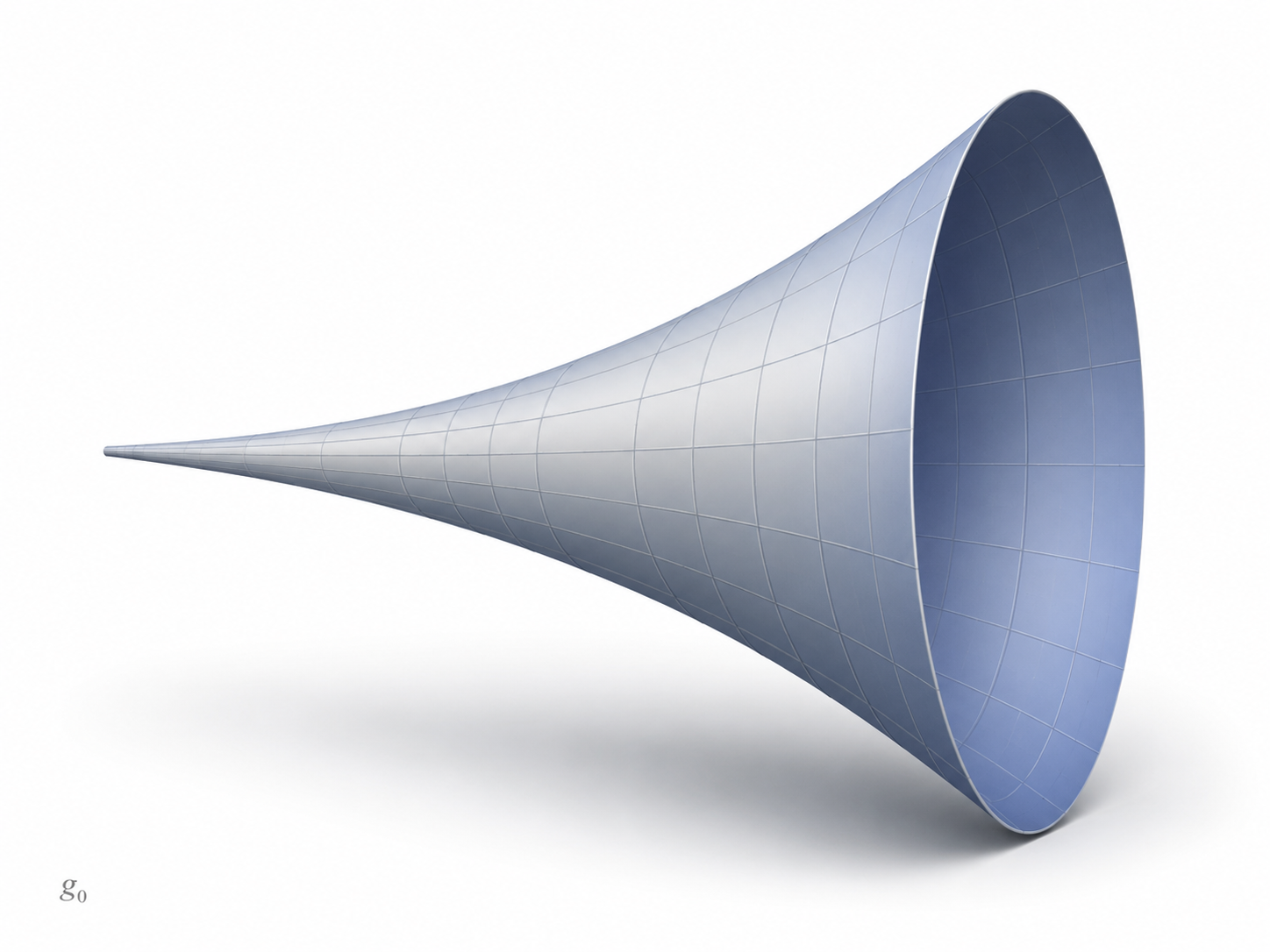}
\end{minipage}\hfill
\begin{minipage}[t]{0.48\textwidth}
\centering
\includegraphics[width=\textwidth]{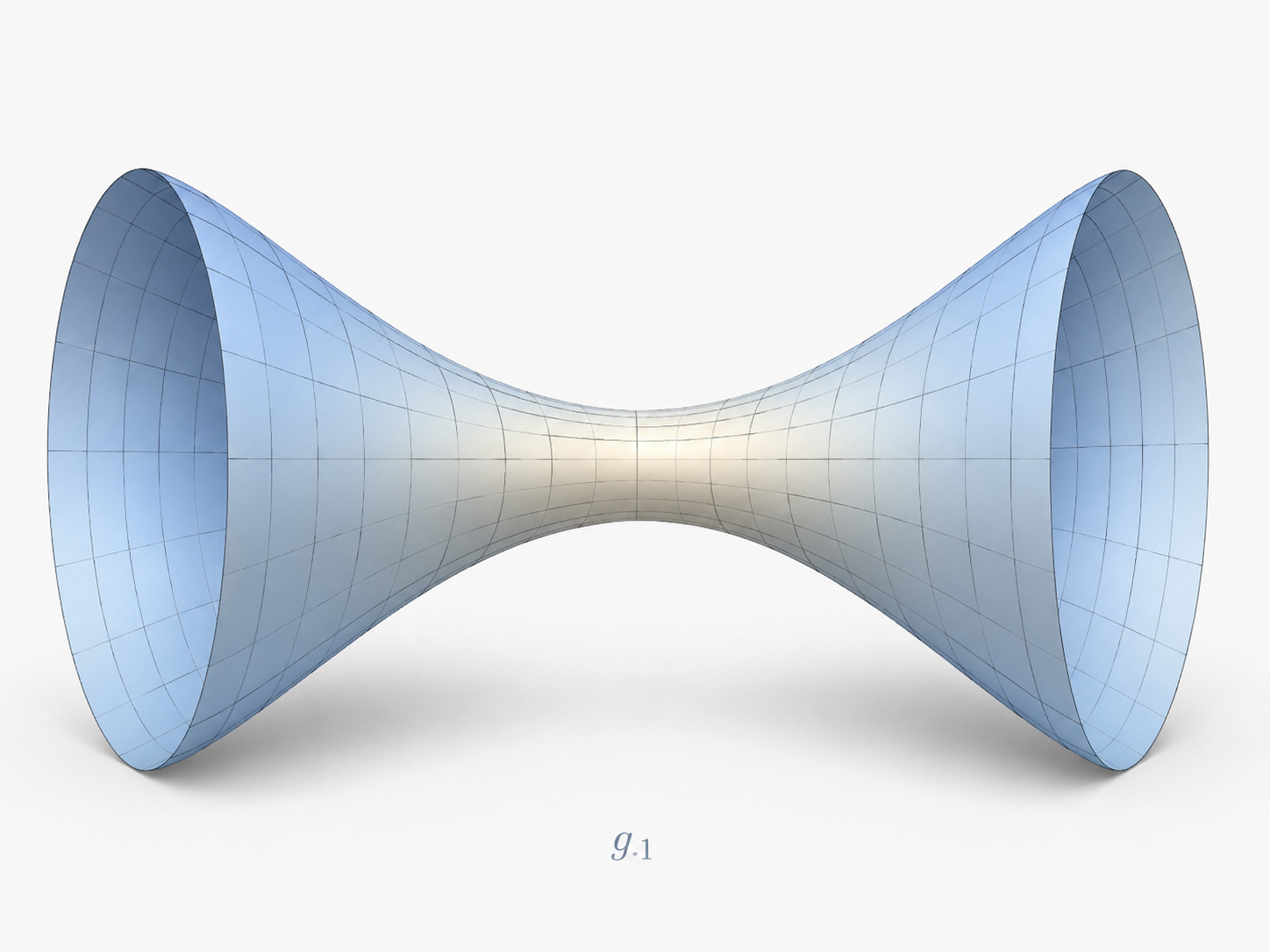}
\end{minipage}
\caption{Schematic renderings of the two warped products. On the left,
$g_0=dx^2+e^{2x}d\theta^2$, with no class $\beta$ geodesic. On the right,
$g_1=dx^2+\cosh^2(x)d\theta^2$, 
with a unique class $\beta$ geodesic at the waist $x=0$.}
\label{fig:surfaces}
\end{figure}

The weak Whitney path described in the introduction is especially
simple in these coordinates. It is given by
\begin{equation*}
f_s(x)=(1-s)e^x+s\cosh x,
\qquad
g_s=dx^2+f_s(x)^2d\theta^2.
\end{equation*}
Since $f_s''=f_s$, every $g_s$ has curvature $-1$.

\section{The cylinder quotient}\label{sec_cylinder_moduli}
\begin{proof}[Proof of Theorem \ref{thm_cylinder_moduli}]
Let $\beta$ be the generator of $\pi_1(\mathbb R\times S^1)$, and put
\begin{equation*}
G=\Diff(\mathbb R\times S^1),\qquad
\mathcal B=\{\beta\}.
\end{equation*}
Every diffeomorphism sends $\beta$ to $\beta$ or $\beta^{-1}$, so
Theorem \ref{thm_mainN} gives a locally constant $\{0,1\}$-valued function
\begin{equation*}
\mathcal N_{\mathcal B,G}:\mathcal M(\mathbb R\times S^1)\to\{0,1\}.
\end{equation*}
Since $\mathcal N_{\mathcal B,G}$ is locally constant, its fibers
are both open and closed.

The metrics $g_0,g_1$ of Section \ref{sec_Acalculation} satisfy
\begin{equation*}
\mathcal N_{\mathcal B,G}([g_0])=0,\qquad
\mathcal N_{\mathcal B,G}([g_1])=1.
\end{equation*}
Therefore $\mathcal M(\mathbb R\times S^1)$ contains two non-empty
disjoint clopen sets, and is disconnected.
\end{proof}

\section{Finite-type surfaces}\label{sec_finite_type}
\begin{proof}[Proof of Theorem \ref{thm_finite_type_moduli}]
Let $p_1,\ldots,p_n$ be the punctures of $X$, and let $\beta_i$ be
the primitive peripheral class going around $p_i$. We first construct,
for each vector
\begin{equation*}
\varepsilon=(\varepsilon_1,\ldots,\varepsilon_n)\in\{0,1\}^n,
\end{equation*}
a metric $g_\varepsilon\in\Neg(X)$ such that
\begin{equation}\label{eq_Fpattern}
N(g_\varepsilon,\beta_i)=\varepsilon_i.
\end{equation}

Choose a pair-of-pants decomposition of a compact core of $X$, whose
external cuffs correspond to the punctures. For the external cuff
corresponding to $p_i$, prescribe length $0$ if $\varepsilon_i=0$, and
length $1$ if $\varepsilon_i=1$; here length $0$ means a cusp. Choose
one positive length for each internal cuff, assigning the same length
to the two copies which are to be glued. Hyperbolic pairs of pants with
these prescribed boundary lengths, allowing length $0$ for cusps,
exist; see for instance ~\cite[Chapter 3]{cite_Buser}. Gluing the
internal cuffs with matching lengths gives a hyperbolic compact core.
At each external cuff of positive length, attach the standard
hyperbolic funnel of the same boundary length. The resulting complete
metric has curvature $-1$ everywhere, hence belongs to
$\Neg(X)$.

If $\varepsilon_i=0$, then the class $\beta_i$ is parabolic for the
hyperbolic metric $g_\varepsilon$, so it has no closed geodesic
representative and
\begin{equation*}
N(g_\varepsilon,\beta_i)=0.
\end{equation*}
If $\varepsilon_i=1$, then $\beta_i$ is represented by the core
geodesic of the attached funnel. This representative is unique. Thus
\begin{equation*}
N(g_\varepsilon,\beta_i)=1.
\end{equation*}
This proves \eqref{eq_Fpattern}.

For the pure diffeomorphism group, set
\begin{equation*}
G_{\mathrm p}=\Diff_{\mathrm p}(X),\qquad
\mathcal B_i=\{\beta_i\}.
\end{equation*}
Each $\mathcal B_i$ is preserved by $G_{\mathrm p}$ up to inversion,
so Theorem \ref{thm_mainN} gives locally constant functions
\begin{equation*}
\mathcal N_{\mathcal B_i,G_{\mathrm p}}:
\mathcal M_{\mathrm p}(X)\to\{0,1\}.
\end{equation*}
Equivalently, the map
\begin{equation*}
\overline\Phi([g])=
\big(\mathcal N_{\mathcal B_1,G_{\mathrm p}}([g]),\ldots,
\mathcal N_{\mathcal B_n,G_{\mathrm p}}([g])\big)
\end{equation*}
is a continuous map
\begin{equation*}
\overline\Phi:\mathcal M_{\mathrm p}(X)\to\{0,1\}^n.
\end{equation*}
The metrics $g_\varepsilon$ realise all binary vectors in
$\{0,1\}^n$. Thus the corresponding fibers of
$\overline\Phi$ give at least $2^n$ non-empty disjoint clopen
subsets. Hence the pure quotient has at least $2^n$ connected
components.

For the full diffeomorphism group, a diffeomorphism may permute the
punctures. Thus the ordered vector need not descend, but its unordered
version does. Set
\begin{equation*}
G=\Diff(X),\qquad
\mathcal B=\{\beta_1,\ldots,\beta_n\}.
\end{equation*}
The set $\mathcal B$ is preserved by $G$ up to permutation and
inversion, and Theorem \ref{thm_mainN} gives
\begin{equation*}
\mathcal N_{\mathcal B,G}:\mathcal M(X)\to
\operatorname{Sym}^{n}(\{0,1\}).
\end{equation*}
In particular the function
\begin{equation*}
\overline N([g])=\#\{i\mid N(g,\beta_i)=1\}
\end{equation*}
is locally constant and well-defined on the full quotient:
\begin{equation*}
\overline N:\mathcal M(X)\to\{0,1,\ldots,n\}.
\end{equation*}
The construction above realises each value $N=0,1,\ldots,n$, and so
the full quotient has at least $n+1$ connected components.
\end{proof}

\end{document}